\documentclass[a4paper,12pt]{amsart}
\usepackage{amssymb}
\usepackage{ifthen}
\usepackage{graphicx}
\usepackage{scrtime}

\newtheorem{thm}{Theorem}
\newtheorem{cor}{Corollary}
\newtheorem{lem}{Lemma}

\newtheorem{rem}{Remark}

\newtheorem{conj}{Conjecture}
\newtheorem{prob}{Problem}
\theoremstyle{definition}
\newtheorem{defn}{Definition}
\newtheorem{example}{Example}

\newcounter {own}
\def\theown {\thesection       .\arabic{own}}

\newenvironment{pf}[1][]{%
 \vskip 3mm
 \noindent
 \ifthenelse{\equal{#1}{}}%
  {{\slshape Proof. }}%
  {{\slshape #1.} }%
 }%
{\qed\bigskip}

\newcounter{alphabet}
\newcounter{tmp}
\newenvironment{Thm}[1][]{\refstepcounter{alphabet}%
\bigskip%
\noindent%
{\bf Theorem \Alph{alphabet}}%
\ifthenelse{\equal{#1}{}}{}{ (#1)}%
{\bf .} \itshape}{\vskip 8pt}

\makeatletter
\newcommand{\Ref}[1]{\@ifundefined{r@#1}{}{\setcounter{tmp}{\ref{#1}}\Alph{tmp}}}
\makeatother

\newenvironment{Lem}[1][]{\refstepcounter{alphabet}%
\bigskip%
\noindent%
{\bf Lemma \Alph{alphabet}}%
{\bf .} \itshape}{\vskip 8pt}

\newenvironment{Conje}[1][]{\refstepcounter{alphabet}%
\bigskip%
\noindent%
{\bf Conjecture \Alph{alphabet}}%
{\bf .} \itshape}{\vskip 8pt}

\newcommand{\es}{{\mathcal S}}

\newcommand{\IC}{{\mathbb C}}
\newcommand{\ID}{{\mathbb D}}

\newcommand{\F}{{\mathcal{F}}}

\newcommand{\D}{{\mathbb D}}

\def\be{\begin{equation}}
\def\ee{\end{equation}}

\newcommand{\bee}{\begin{enumerate}}
\newcommand{\eee}{\end{enumerate}}

\newcommand{\blem}{\begin{lem}}
\newcommand{\elem}{\end{lem}}
\newcommand{\bthm}{\begin{thm}}
\newcommand{\ethm}{\end{thm}}
\newcommand{\bcor}{\begin{cor}}
\newcommand{\ecor}{\end{cor}}
\newcommand{\beg}{\begin{example}}
\newcommand{\eeg}{\end{example}}
\newcommand{\begs}{\begin{examples}}
\newcommand{\eegs}{\end{examples}}
\newcommand{\bdefe}{\begin{defn}}
\newcommand{\edefe}{\end{defn}}
\newcommand{\bprob}{\begin{prob}}
\newcommand{\eprob}{\end{prob}}
\newcommand{\bei}{\begin{itemize}}
\newcommand{\eei}{\end{itemize}}

\newcommand{\bcon}{\begin{conj}}
\newcommand{\econ}{\end{conj}}
\newcommand{\bcons}{\begin{conjs}}
\newcommand{\econs}{\end{conjs}}
\newcommand{\bprop}{\begin{propo}}
\newcommand{\eprop}{\end{propo}}
\newcommand{\br}{\begin{rem}}
\newcommand{\er}{\end{rem}}
\newcommand{\brs}{\begin{rems}}
\newcommand{\ers}{\end{rems}}
\newcommand{\bo}{\begin{obser}}
\newcommand{\eo}{\end{obser}}
\newcommand{\bos}{\begin{obsers}}
\newcommand{\eos}{\end{obsers}}
\newcommand{\bpf}{\begin{pf}}
\newcommand{\epf}{\end{pf}}
\newcommand{\ba}{\begin{array}}
\newcommand{\ea}{\end{array}}
\newcommand{\beq}{\begin{eqnarray}}
\newcommand{\beqq}{\begin{eqnarray*}}
\newcommand{\eeq}{\end{eqnarray}}
\newcommand{\eeqq}{\end{eqnarray*}}

\newcommand{\ds}{\displaystyle}

\newcounter{minutes}
\divide\time by 60
\newcounter{hours}
\multiply\time by 60 \addtocounter{minutes}{-\time}

\begin{document}

\bibliographystyle{amsplain}

%

\title[Uniform Close-to-convexity radius of sections of functions in the close-to-convex
family]{Uniform Close-to-convexity radius of sections of functions in the close-to-convex
family}

\def\thefootnote{}
\footnotetext{ \texttt{\tiny File:~\jobname .tex,
          printed: \number\day-\number\month-\number\year,
          \thehours.\ifnum\theminutes<10{0}\fi\theminutes}
} \makeatletter\def\thefootnote{\@arabic\c@footnote}\makeatother

\author{Vaidhyanathan Bharanedhar}
\address{S. V. Bharanedhar, Indian Statistical Institute (ISI),
 Chennai Centre, SETS (Society
for Electronic Transactions and security), MGR Knowledge City, CIT
Campus, Taramani, Chennai 600 113, India.}
\email{bharanedhar3@gmail.com}

\author{Saminathan Ponnusamy$^\dagger $
}
\address{S. Ponnusamy,
Indian Statistical Institute (ISI), Chennai Centre, SETS (Society
for Electronic Transactions and security), MGR Knowledge City, CIT
Campus, Taramani, Chennai 600 113, India. }
\email{samy@isichennai.res.in, samy@iitm.ac.in}


\keywords{Univalent, convex and close-to-convex functions;
 partial sums, convolution.  \\
$
^\dagger$ {\tt The first author is currently on leave from the
Department of Mathematics, Indian Institute of Technology Madras, Chennai-600 036, India.}
}

\begin{abstract}
The authors consider the class $\F$ of normalized functions $f$ analytic in the unit disk $\ID$
and satisfying the condition
$${\rm Re}\left(1+\frac{zf''(z)}{f'(z)}\right)>-\frac{1}{2},\quad z\in\D.
$$
Recently, Ponnusamy et al. \cite{samy-hiroshi-swadesh} have shown that $1/6$ is the uniform sharp bound for the
radius of convexity of every section of each function in the class $\F$. They conjectured that $1/3$ is the
uniform univalence radius of every section of $f\in \F$. In this paper, we solve this conjecture affirmatively.

\end{abstract}

\maketitle
\pagestyle{myheadings}
\markboth{S. V. Bharanedhar and S. Ponnusamy}
{Sections of functions in the close-to-convex
family}

\section{Preliminaries and the Main Theorem}
Let ${\mathcal A}$ be the family of functions analytic in the unit disk
$\ID := \{z\in \IC:\, |z| <1 \}$ of the form
$f(z)=z+\sum_{k=2}^\infty a_k z^k.
$
Then the $n$-th \emph{section/partial sum} of $f$,
denoted by $s_n(f)(z)$, is  defined to be the polynomial
$$s_n(f)(z)=z+\sum_{k=2}^n a_k z^k.
$$ 
Let ${\mathcal S}$ denote the class of functions  in ${\mathcal A}$ that are  univalent in $\ID $.
Finally, let $\mathcal{C}$, $\mathcal{S}^*$ and $\mathcal{K}$  denote the usual geometric subclasses of functions in
$\es$ with convex, starlike and close-to-convex images, respectively (see \cite{Duren:univ}).

If $f\in {\mathcal S}$ is arbitrary, then the argument principle shows that the
$n$-th section $s_n(f)(z)$ is univalent in each fixed compact disk $|z|\leq r~(<1)$ provided that
$n$ is sufficiently large. But then if we set $p_n(z)=r^{-1}s_n(f)(rz)$, then $p_n(z)$ is a
polynomial that is univalent in the unit disk $\ID$. Consequently, the set of
univalent polynomials is dense with respect to the topology of locally uniformly
in $\mathcal S$ (see \cite{Duren:univ}. Suffridge \cite{Suff72} showed that even the subclass of
polynomials with the highest coefficient $a_n=1/n$ is dense in $\mathcal S$.
Szeg\"o \cite{szego28} discovered that every section $s_n(f)$
is univalent in the disk $|z|<1/4$ for all  $f\in\es$ and for each $n\geq2$.
The radius $1/4$ is best possible as the Koebe function $k(z)=z/(1-z)^2$ shows.
It is worth pointing out that the case $n=3$ of Szeg\"o's result is far from triviality.

In \cite{Rush88}, Ruscheweyh established a stronger result by showing that the partial sums
$s_n(f)(z)$ of $f$ are indeed starlike in the disk $|z|<1/4$ for functions $f$ belonging not
only to  $\mathcal{S}$ but also to the closed convex hull of $\mathcal{S}$.
The following conjecture concerning the exact (largest) radius of univalence $r_n$ of $f\in {\mathcal S}$
is still open (see \cite{Robertson} and \cite[\S8.2, p.~241--246]{Duren:univ}).

\begin{Conje}
If $f\in\es$, then $s_n(f)$ is univalent in $|z|<1-\frac{3}{n}\log n$
for all $n\geq5$.
\end{Conje}

A surprising fact observed by Bshouty and Hengartner \cite{BH1991}
is that the Koebe function is no more extremal for the above conjecture.
On the other hand, this conjecture has been solved by using an important
convolution theorem \cite{Rush73} for a number of geometric subclasses  of $\es$, for example,
the classes $\mathcal{C}$, $\mathcal{S}^*$ and $\mathcal{K}$.
Indeed, for $\phi (z)=z/(1-z)$, the sections $s_n(\phi)$ are known to be convex in $|z|<1/4$ (see \cite{GoodScho}).
Moreover for the Koebe function  $k(z)=z/(1-z)^2$,  $s_n(k)$ is known to be starlike
in $|z|<1-\frac{3}{n}\log n$ for $n\geq 5$ and hence, for the convex function $\phi (z)=z/(1-z)$,
$s_n(\phi)$ is convex in $|z|<1-\frac{3}{n}\log n$ for $n\geq 5$.
From a convolution theorem relating to the P\'{o}lya-Schoenberg conjecture proved by Ruscheweyh and  Sheil-Small
\cite{Rush73}, it follows that all sections $s_n(f)$ are  convex (resp. starlike, close-to-convex) in $|z|<1/4$ whenever
$f\in\mathcal{C}$ (resp. $f\in \mathcal{S}^*$ and $f\in \mathcal{K}$).  Similarly, for  $n\geq 5$,
$s_n(f)$  is  convex (resp. starlike, close-to-convex) in $|z|<1-\frac{3}{n}\log n$ whenever
$f\in\mathcal{C}$ (resp. $f\in \mathcal{S}^*$ and $f\in \mathcal{K}$).
An account of history of this and related information may be found in \cite[\S8.2, p.~241--246]{Duren:univ}
and also in the nice survey article of Iliev  \cite{iliev}.   For further interest on this topic, we refer to
\cite{FourSilver91,Rus-72, Small-70,Silver88} and recent articles \cite{OP11-SMJ,MP12-preprint,MP12-RMJM,OPW13-SMJ}.


%
One of the important criteria for an analytic function $f$ defined on a convex domain
$\Omega$, to be univalent in $\Omega$ is that ${\rm Re}\,f'(z)>0$ on $\Omega$
(see \cite[Theorem 2.16, p. 47]{Duren:univ}). The following definition is a consequence of it.

A function $f\in\mathcal{A}$ is said to be close-to-convex (with respect to $g$),
denoted by $f\in\mathcal{K}_g$ if  there exists a $g\in\mathcal{C}$
such that
\begin{equation}\label{Bpon6-eq1c}
{\rm Re}\left(e^{i\alpha}\frac{f'(z)}{g'(z)}\right)>0,\,z\in\ID,
 \end{equation}
for some real $\alpha$ with  $|\alpha|<\pi/2$. More often, we
consider $\mathcal{K}_g$ (with $\alpha=0$ in \eqref{Bpon6-eq1c})
 and
 $\mathcal{K}=\mathop{\cup}_{g\in\mathcal{C}}\mathcal{K}_g. $
For functions
in $\mathcal{K}_g$, we have the following result of Miki \cite{Miki}.

\begin{Thm}\label{BPO2Tha2}
Let $f\in\mathcal{K}_g$, where $g(z)=z+\sum_{n=2}^{\infty}b_nz^n$. Then $s_n(f)$ is close-to-convex
with respect to $s_n(g)$ in $|z|<1/4$.
\end{Thm}

In a recent paper \cite{bhaobpo}, the present authors proved the following.

\begin{Thm}\label{Bpon6-thbhaobpo}
Let $f\in\mathcal{K}$. Then every section $s_{n}(f)$ of $f$ belongs
to the class $\mathcal{K}$  in the disk
$|z|<1/2$ for all $n \geq 46$.
\end{Thm}
Choosing different convex functions $g$
in \cite{bhaobpo}, the authors have found the
value $N(g)\in\mathbb{N}$ for $f\in\mathcal{K}_g$ such that
$s_n(f)\in\mathcal{K}_g$ in a disk $|z|<r$ for all $n\geq N(g)$.

In \cite{samy-hiroshi-swadesh}, the authors consider the
class $\F$ of locally univalent functions $f$ in $\mathcal A$
satisfying the condition
\begin{equation}\label{Bpon6-eq1b}
{\rm Re}\left(1+\frac{zf''(z)}{f'(z)}\right)>-\frac{1}{2},\quad z\in\D.
\end{equation}

The importance of this class is outlined in \cite{samy-hiroshi-swadesh} and it
was also remarked  that the class $\F$ has a special
role on certain problems on the class of harmonic univalent mappings in $\ID$ (see \cite{samy-hiroshi-swadesh} and
the references therein). It is worth remarking that functions in $\F$ are neither included in $\es^*$
nor includes $\es^*$ nor $\mathcal{K}$.
It is well-known that $\F\subsetneq {\mathcal K}\subsetneq {\mathcal S}$
and hence, it is obvious from an earlier observation that for $f\in\F$, each $s_n(f)(z)$  is close-to-convex in $|z|<1/4$.
An interesting question is to determine the largest uniform disk with this property (see Conjecture \ref{conj1} below).
We now recall a recent result of Ponnusamy et al. \cite{samy-hiroshi-swadesh}.

 \begin{Thm}\label{theorem-for-convexity}
 Every section of a function in the class $\F$ is convex in the disk $|z|<1/6$.
  The radius $1/6$
 cannot be replaced by a greater one.
 \end{Thm}

In the same article the authors \cite{samy-hiroshi-swadesh}
observed that all sections functions of $\F$ are close-to-convex in
the disk $|z|<1-\frac{3}{n}\log n$ for $n\ge 5$.  Consider
\begin{equation}\label{Bpon6-eq1a}
f_0(z)=\frac{z-z^2/2}{(1-z)^2}.
\end{equation}
We see that $f_0\notin{\mathcal S}^*$, but $f_0\in\mathcal{K}$.
Also,  $f_0$ is extremal for many extremal problems for the class $\mathcal{F}$.
By investigating the second partial sum of $f_0\in\F$, the authors conjectured the following.

\begin{conj}\label{conj1}
Every section $s_n(f)$ of $f\in \F$ is close-to-convex in the disk $|z|<1/3$ and $1/3$ is sharp.
\end{conj}

In this article we solve this conjecture in  the following form.

\begin{thm}\label{thm1}
Every section $s_n(f)$ of $f\in \F$ satisfies ${\rm Re}\,(s_n(f)'(z))>0$
in the disk $|z|<1/3$. In particular every section is close-to-convex in the disk $|z|<1/3$.
The radius $1/3$ cannot be replaced by a greater one.
\end{thm}


We remark that this result is much stronger than the original conjecture. The following lemma is useful in the proof of
Theorem~\ref{thm1}.

\begin{Lem}\cite[Lemma 1]{samy-hiroshi-swadesh}\label{lemma-for-convexity}
If $f(z)=z+\sum_{n=2}^\infty a_nz^n\in \F$, then the following estimates
hold:
\begin{enumerate}
\item[{\bf (a)}] $|a_n|\le \ds\frac{n+1}{2}$ for  $n\ge 2$. Equality holds for
$ f_0(z)$ given by \eqref{Bpon6-eq1a} or its rotation.

\item[{\bf (b)}] $\ds\frac{1}{(1+r)^3}\le |f'(z)|\le \frac{1}{(1-r)^3}$ for  $|z|=r<1$.
The bounds are sharp.

\item[{\bf (c)}] If $f(z)=s_n(z)+\sigma_n(z)$,
with $\sigma_n(z)=\sum_{k=n+1}^\infty a_kz^k$, then for $|z|=r<1$
we have
$$\ds|\sigma_n'(z)|\le \frac{n(n+1)r^{n+2}-2n(n+2)r^{n+1}+(n+1)(n+2)r^n}{2(1-r)^3}.
$$
\end{enumerate}
\end{Lem}

\section{Proof of Theorem~\ref{thm1}}
Let $f(z)=z+\sum_{n=2}^\infty a_nz^n\in \F$.
We shall prove that each partial sum $s_n(z):=s_n(f)(z)$ of $f$
satisfies the condition ${\rm Re}\,(s_n'(z))>0$  in the disk $|z|<1/3$
for all $n\geq2$.

Let us first consider the second section $s_2(z)=z+a_2z^2$ of $f$.
A simple computation shows that
$${\rm Re}\,(s_2'(z))=1+{\rm Re}\,(2a_2z).
$$
From Lemma~\Ref{lemma-for-convexity}(a),
we have $|a_2|\leq 3/2$ and as a consequence of it we get
$${\rm Re}\,(s_2'(z))\geq 1-2|a_2||z|\geq1-3|z|
$$
which is positive provided $|z|<1/3$. Thus, $s_2(z)$ is
close-to-convex in the disk $|z|<1/3$.
To show that the constant $1/3$ is best possible, we consider the function
$f_0\in\F$ given in \eqref{Bpon6-eq1a}, namely,
$$f_0(z)=\frac{1}{2}\left[\frac{1}{(1-z)^2}-1\right] =
z+\sum_{n=2}^{\infty}\left (\frac{n+1}{2}\right )z^n.
$$
Let us denote by $s_{2,0}(z)$,
the second partial sum $s_2(f_0)(z)$ of $f_0(z)$
so that $s_{2,0}(z)=z+(3/2)z^2$. Then we get
$s_{2,0}'(z)=1+3z,
$
which vanishes at  $z=-1/3$. Thus the constant $1/3$ is best possible.

Next, let us consider the case $n=3$.
Each $f\in\F$ satisfies the analytic condition \eqref{Bpon6-eq1b} and so we can write
\begin{equation}\label{thm1eq0}
1+\frac{2}{3}\frac{zf''(z)}{f'(z)}=p(z),
\end{equation}
where $p(z)=1+p_1z+p_2z^2+\cdots$ is analytic in
 $\D$ and ${\rm Re}\, p(z)>0$ in $\D$. From
Carath\'eodory Lemma~\cite[p.~41]{Duren:univ}
we get $|p_n|\le 2$ for all $n\ge 2$.
If we rewrite (\ref{thm1eq0}) in power series form, then
$$1+\frac{2}{3}\frac{z(2a_2+6a_3z+12a_4z^2+\cdots)}{1+2a_2z+3a_3z^2+\cdots}
=1+p_1z+p_2z^2+\cdots.
$$
Now comparing the coefficients of $z$ and $z^2$ on both sides yields the relations
$$p_1=\frac{4}{3}a_2 ~~\mbox{ and }~~ p_2=\frac{4}{3}(3a_3-2a_2^2).
$$
As $|p_1|\le 2$ and $|p_2|\le 2$, we may rewrite the last two relations as
\begin{equation}\label{thm1eq-0}
a_2=\frac{3}{2}\alpha ~~\mbox{ and }~~ \frac{2}{3}(3a_3-2a_2^2)=\beta,
~\mbox{i.e. $a_3=\displaystyle\frac{1}{2}(\beta+3\alpha^2)$}
\end{equation}
for some $|\alpha|\le 1$ and $|\beta|\le 1$. Now we have to show that
\begin{equation}\label{Bpon6-eq2}
{\rm Re}\,(s_3'(z))={\rm Re}\,(1+2a_2z+3a_3z^2)>0
\end{equation}
in $|z|<1/3$. Since the function ${\rm Re}\,(s_3'(z))$
is harmonic in $|z|\leq1/3$, it is enough to prove \eqref{Bpon6-eq2}
for $|z|=1/3$. By considering a suitable rotation of $f$, it is enough to prove
\eqref{Bpon6-eq2} for $z=1/3$. Thus, it suffices to show that
\begin{equation}\label{Bpon6-eq3}
{\rm Re}\left(1+\frac{2}{3}a_2+\frac{1}{3}a_3\right)>0.
\end{equation}
By using the relations in \eqref{thm1eq-0} and the maximum principle, we see that
the inequality \eqref{Bpon6-eq3} is equivalent to
\begin{equation}\label{Bpon6-eq4}
{\rm Re}\left(1+\alpha+\frac{\alpha^2}{2}+\frac{\beta}{6}\right)>0,
\end{equation}
where $|\alpha|=1$ and $|\beta|=1$. If we take $\alpha=e^{i\theta}$
and $\beta=e^{i\phi}$ $(0\leq\theta,\phi<2\pi)$, then in order to verify
the inequality  \eqref{Bpon6-eq4} it suffices to prove
$$\displaystyle \min_{\theta,\phi}T(\theta,\phi)>0,
$$
where
$$T(\theta,\phi)=1+\cos\theta+\frac{\cos2\theta}{2}
+\frac{\cos\phi}{6}
$$
and $\theta$, $\phi$ lies in $[0, 2\pi)$. Let
$$g(\theta)=1+\cos\theta+\frac{\cos2\theta}{2}, \quad \theta \in [0, 2\pi).
$$
Then
$$g'(\theta)=-\sin\theta(1+2\cos\theta)~\mbox{ and }~ g''(\theta)=-[\cos\theta+2\cos2\theta].
$$
The points  at which $g'(\theta)=0$ are $\theta=0,\,2\pi/3,\, \pi$ and $4\pi/3$.
But $g''(\theta)$ is positive for $\theta=2\pi/3$ and $\theta=4\pi/3$.  Hence
$$\min_{\theta}g(\theta)=g\left(\frac{2\pi}{3}\right)=g\left(\frac{4\pi}{3}\right)=\frac{1}{4}.
$$
As the minimum value of $(\cos\phi)/6$ is $-1/6$, it follows that
$$\min_{\theta,\phi}T(\theta,\phi)=T\left(\frac{2\pi}{3},\pi\right)
=T\left(\frac{4\pi}{3},\pi\right)=\frac{1}{12}>0.
$$
This proves the inequality \eqref{Bpon6-eq2} for $|z|<1/3$.

Now let us consider the case $n\geq4$. Let $f(z)=s_n(z)+\sigma_n(z)$,
where $\sigma_n(z)$ is as given in Lemma \Ref{lemma-for-convexity}(c).
Then 
\begin{equation}\label{Bpon6-eq4a}
{\rm Re}\,(s_n'(z))={\rm Re}\,(f'(z)-\sigma_n'(z))\geq {\rm Re}\,(f'(z))-|\sigma_n'(z)|.
\end{equation}
By maximum principle it is enough to prove that ${\rm Re}\,(s_n'(z))>0$
for $|z|=1/3$. Now let us estimate the values of ${\rm Re}\,(f'(z))$ and
$|\sigma_n'(z)|$ on $|z|=1/3$.

\begin{figure}
\begin{center}
\includegraphics[height=6cm, width=6cm, scale=1]{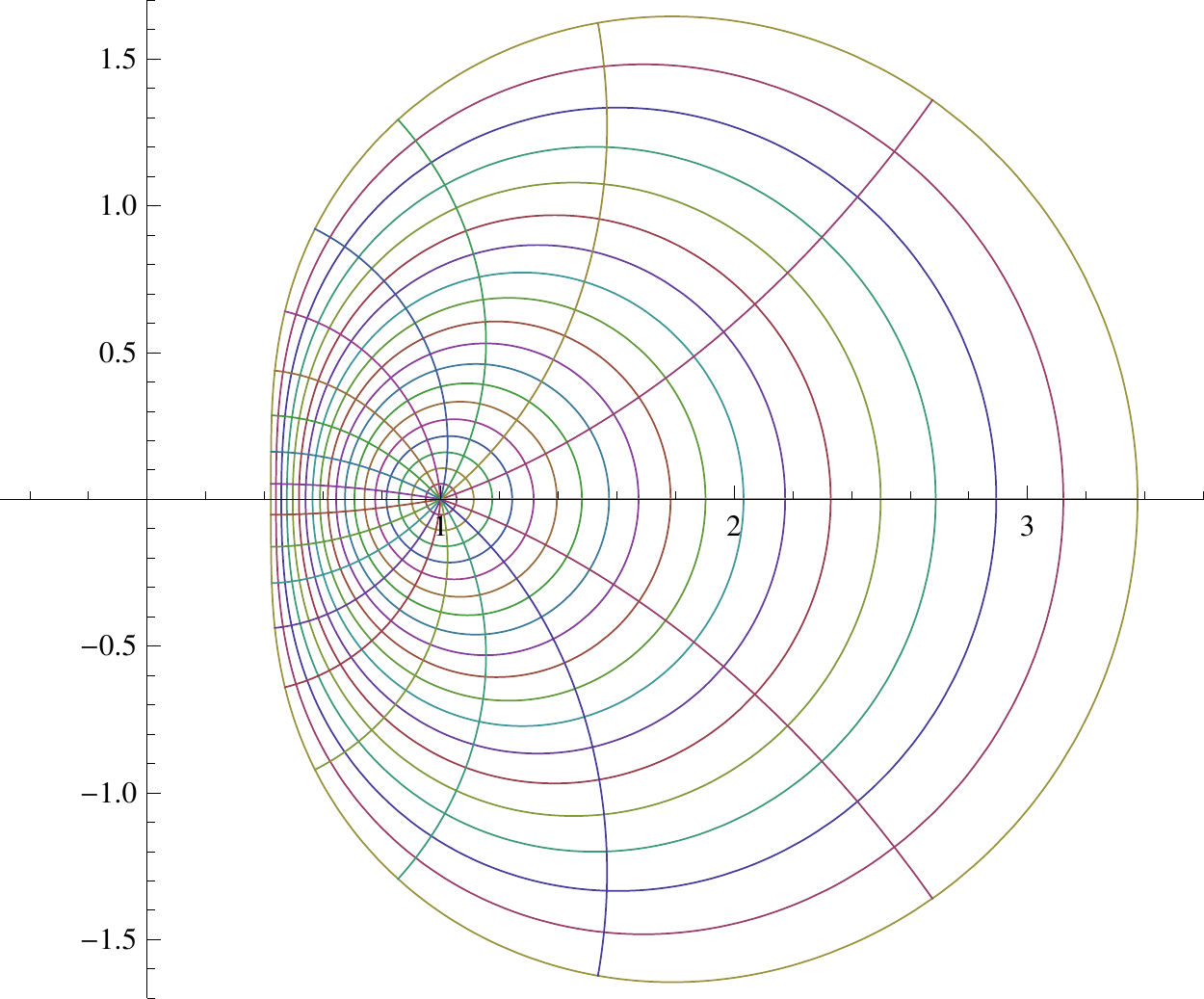}
\hspace{1cm}
\includegraphics[height=6cm, width=6cm, scale=1]{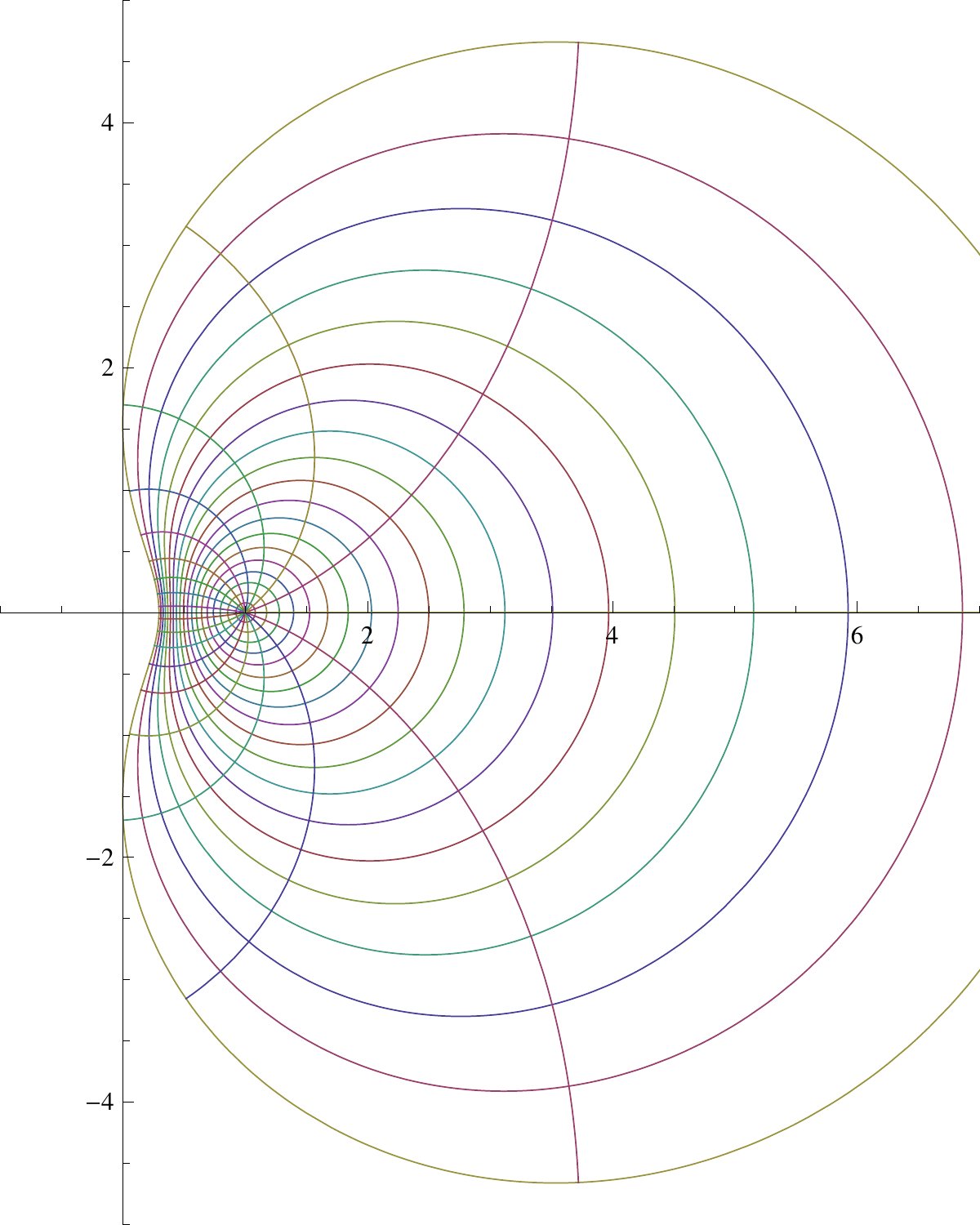}
\end{center}
\hspace{0.4cm} \textbf{(a)}: $H(|z|<1/3)$ \hspace{4.5cm} \textbf{(b)}: $H(|z|<1/2)$
\vspace{2 cm}
\begin{center}
\includegraphics[height=6cm, width=6cm, scale=1]{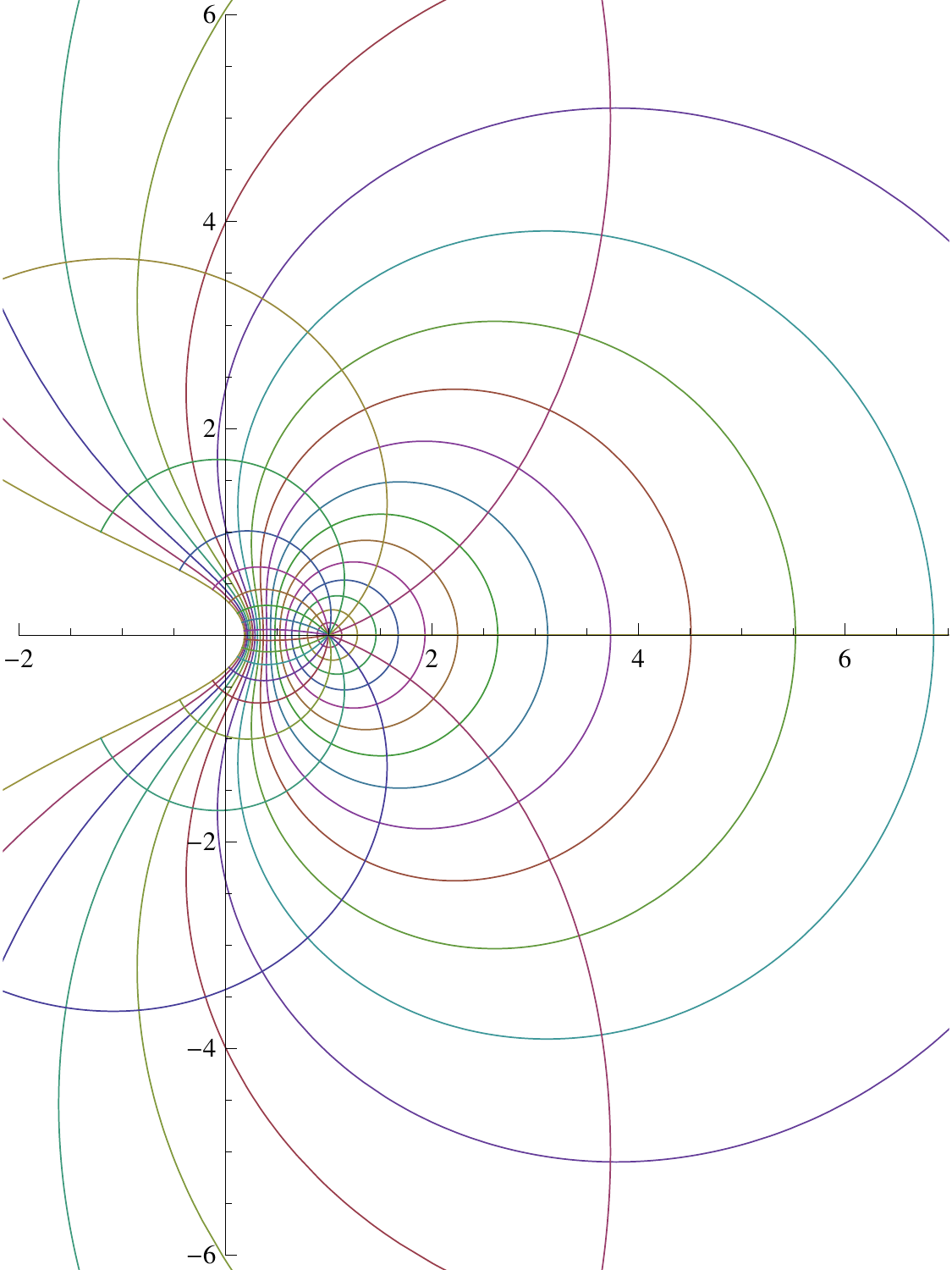}
\hspace{1cm}
\includegraphics[height=6cm, width=6cm, scale=1]{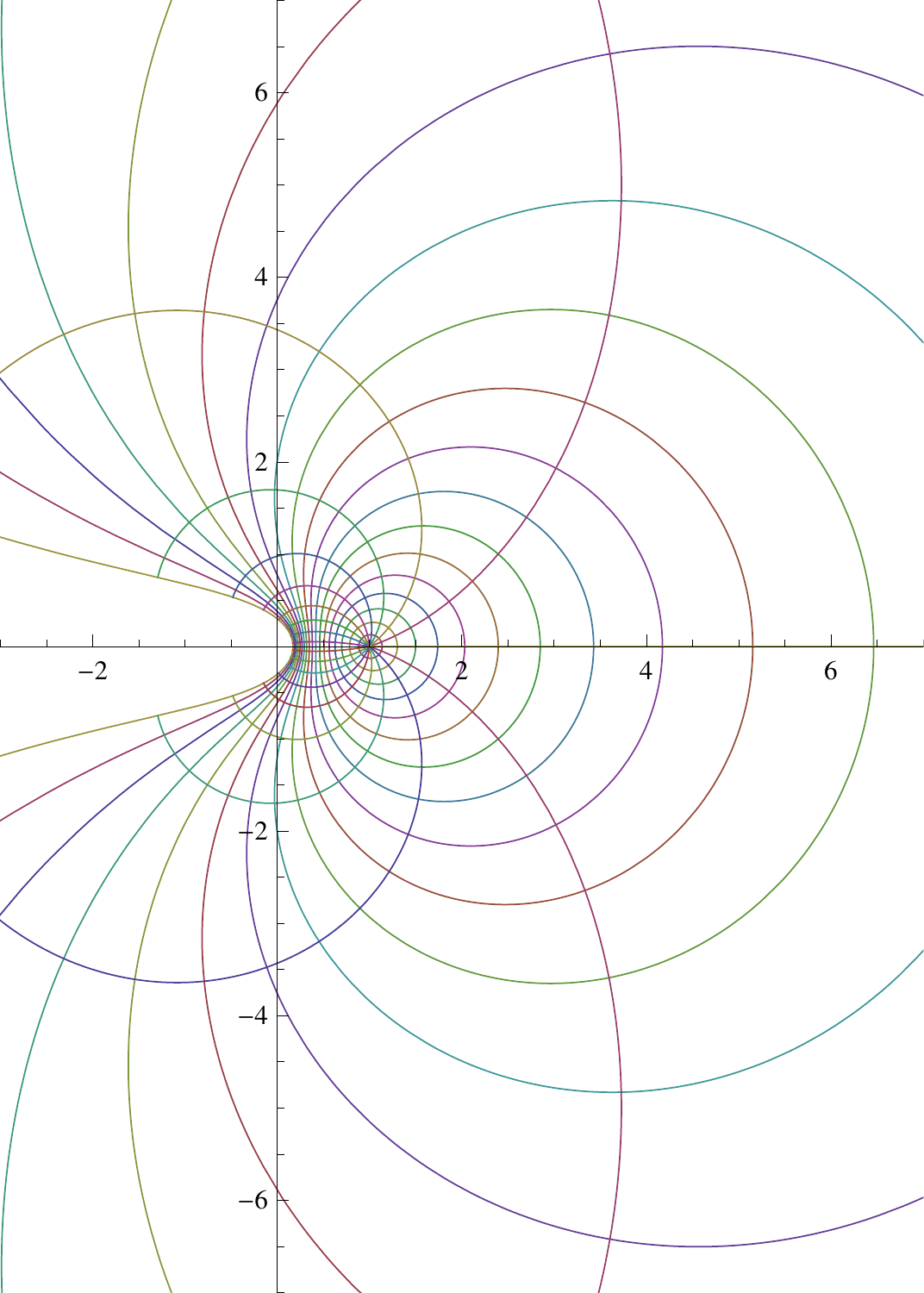}
\end{center}
\hspace{0.4cm}\textbf{(c)}: $H(|z|<3/4)$ \hspace{4.5cm} \textbf{(d)}: $H(|z|<4/5)$
\caption{Images of the disks $|z|<r$ for $r=1/3,1/2,3/4,4/5$, under the
function $H(z)=1/(1-z)^3$ \label{subordination}.}
\end{figure}

As in the proof of Lemma \Ref{lemma-for-convexity}(b) in
\cite{samy-hiroshi-swadesh}, we have the subordination relation
for $f\in\F$,
\begin{equation}\label{Bpon6-eq5}
f'(z)\prec\frac{1}{(1-z)^3},~z\in\ID.
\end{equation}
We need to find the image of the circle $|z|=r$ under the transformation
$\ds w(z)=1/(1-z)^3$. As the bilinear transformation $T(z)=1/(1-z)$ maps the circle $|z|=r$ onto the circle
$$\left|T-\frac{1}{1-r^2}\right|=\frac{r}{1-r^2},~\mbox{i.e.,}~T(z)=
\frac{1+re^{i\theta}}{1-r^2},
$$
a little computation shows that the image of the circle $|z|=r$ under
the transformation $\ds w= 1/(1-z)^3$ is a closed curve described by
$$w=\frac{(1+re^{i\theta})^3}{(1-r^2)^3}
=\frac{1+r^3e^{3i\theta}+3r^2e^{2i\theta}+3re^{i\theta}}{(1-r^2)^3},  \quad \theta \in [0, 2\pi) .
$$
From this relation, the substitution  $r=1/3$ gives that
$${\rm Re}\,w=\left(\frac{9}{8}\right)^3
\left[1+\cos\theta+\frac{\cos2\theta}{3}+\frac{\cos3\theta}{27}\right]=h(\theta)~\mbox{(say)}.
$$
If we write $h(\theta)$ in powers of $\cos\theta$, then we easily get
$$h(\theta)=\left(\frac{9}{8}\right)^3\left[ \frac{2}{3}+\frac{8}{9}\cos\theta
+\frac{2}{3}\cos^2\theta+\frac{4}{27}\cos^3\theta \right].
$$
If we let $x=\cos\theta$, then we can rewrite $h(\theta)$  in terms of $x$ as
$$p(x)=\left(\frac{9}{8}\right)^3\left[ \frac{2}{3}+\frac{8}{9}x
+\frac{2}{3}x^2+\frac{4}{27}x^3 \right],
$$
where $-1\leq x\leq 1$. In order to find the minimum value of $h(\theta)$ for $\theta\in[0, 2\pi)$,
it is enough to find the minimum value of $p(x)$ for $x\in[-1,1]$. A computation shows that
$$p'(x)=\frac{81(x+2)(x+1)}{128}
~\mbox{ and }~ p''(x)=\frac{81(3+2x)}{128}.
$$
In the interval $[-1,1]$, $p'(x)=0$ implies $x=-1$ is the only possibility. Also $p''(-1)>0$ and so the minimum value
of the function $p(x)$ in $[-1, 1]$ occurs at $x=-1$. The above discussion implies that
$$\min_{\theta}h(\theta)=h(\pi)=\frac{27}{64}.
$$
Moreover, from the subordination relation \eqref{Bpon6-eq5}, we deduce that
\begin{equation}\label{Bpon6-eq6}
\min_{|z|=1/3}{\rm Re}\,(f'(z))\geq\min_{|z|=1/3}{\rm Re}\left(\frac{1}{(1-z)^3}\right)
=\frac{27}{64}.
\end{equation}
Images of the disks $|z|<r$ for $r=1/3,1/2,3/4,4/5$, under the
function $H(z)=1/(1-z)^3$ are drawn in Figures \ref{subordination}\textbf{(a)}-\textbf{(d)}.
From Lemma \Ref{lemma-for-convexity}(c), we have for $|z|=1/3$
\begin{equation}\label{Bpon6-eq7}
-|\sigma_n(z)|\geq\frac{-1}{8\times 3^{n-1}}
\left[2n^2+8n+9\right]=k(n)~\mbox{(say)}.
\end{equation}
Now
$$k'(n)=\frac{-1}{8\times 3^{n-1}}\left[\log\left(\frac{1}{3}\right)
\left(2n^2+8n+9\right)+4n+8\right].
$$
For $n\geq4$, $k'(n)>0$ and hence $k(n)$ is an increasing function of $n$. Thus for all $n\geq4$,
we have $k(n)\geq k(4)=-73/216$.

Finally, from the relations \eqref{Bpon6-eq4a}, \eqref{Bpon6-eq6} and \eqref{Bpon6-eq7}
it follows that
$${\rm Re}\,(s_n'(z))>\frac{27}{64}-\frac{73}{216}=\frac{145}{1728}>0~\mbox{for all }n\geq4.
$$
The proof is complete.

\vspace{2 cm}
We end the paper with the following conjecture.
\begin{conj}\label{conj2}
Every section $s_n(f)$ of $f\in \F$ is starlike in the disk $|z|<1/3$.
\end{conj}

\end{document}